\documentclass[12pt]{article}
\usepackage{graphicx}
\usepackage{amsmath,amsfonts}
\usepackage{amssymb}

\title{Estimated optimality and robustness of nonlinear  adaptive control systems under bounded disturbances
\author{Alexander Fradkov\\
 Institute for Problems of Mechanical Engineering of RAS,\\
  Saint Petersburg,\\
 Saint Petersburg University\\
 E-mail: fradkov@mail.ru}
}
\date{March 2025}

\newcommand{\trn}{^{\rm\scriptscriptstyle T}}
\renewcommand{\kappa}{\varkappa}
\renewcommand{\le}{\leqslant}
\renewcommand{\ge}{\geqslant}
\renewcommand\epsilon{\varepsilon}

\begin{document}
\maketitle

\begin{abstract}
The problem of suboptimality under bounded disturbances for the adaptive systems based on speed-graadient approach is discussed. A formulation of the estimated optimality of nonlinear nonlinearly parametrized adaptive control systems
is given and sufficient conditions for the estimated optimality in a specified uncertainty class are described for algorithms robustified by negative parametric feedback. A special case of passification based adaptive control for linear time-invariant systems is studied based on the speed-gradient methodology.
\end{abstract}



\section{Introduction}
\par   The evaluation of the adaptive systems
performance under disturbances is a long standing problem since late 1970s \cite{AnnaFrad21}.
Already in the 1970 some ways to achieve robustness under distrubances were proposed. E.g.
inserting negative parametric feedback into the adaptation algorithm (also known as ``$\sigma$-modification'')
was shown to be an efficient way to achieve robustness under bounded input disturbances
both for linear systems \cite{Narendra71,Lindorff73} and for nonlinear systems \cite{F80} .
Improvement of the adaptive systems robustness attracted new attention recently \cite{Karafyllis24,Karafyllis25,KarafyllisKrstic24}.

A more difficult problem is ensuring a suboptimal performance under disturbances, i.e providing
performance arbitrarily close to that of an optimal system. Some results in this direction
were presented in the book \cite{Fradkov90}. However the book \cite{Fradkov90} is written in Russian
and a number of its results were never published in English. The goal of this note is to introduce the
the results on suboptimal (or estimated suboptimal) adaptive control to the worldwide community.

A general result is formulated for a class of nonlinear nonlinearly parametrized problems
where a nonlinear nonlinearly parametrized plant is controlled by the speed-gradient algorithms
(Theorem 1 in Section 3). Then the case of adaptive control of a linear passifiable plant
by adaptive output feedback controller is studied as example.

This paper is dealing with the framework introduced in \cite{F80} where the so called {\it speed-gradient (SG) algorithms} for nonlinear and adaptive control were proposed. It was extended and employed to numerous theoretical and practical problems, see e.g. recent survey \cite{AF21} and more than 200 references therein. Therefore we start with
a brief exposition of the idea of the speed-gradient method.

\section{Speed-gradient method}

Let us recall  briefly the conventional speed-gradient design scheme \cite{F80,AF21}. Consider a nonlinear controlled system in state space described by a system of differential equations
\begin{align}
\dot x(t) = F\big(x, \theta,t \big),
\label{sg1}
\end{align}
where  $x,\in R^n$ is the state vector of the controlled system; $\theta \in R^N$ is the vectors of inputs (controlling or adjustable variables). Let the control goal be specified as follows
\begin{align}
\lim_{t\to\infty}Q(x(t)) = 0 ,
\label{sg2}
\end{align}
where $Q(x)\ge 0$ is a nonnegative smooth goal function.

The speed-gradient approach is aimed at finding the adaptation law for $\theta(t)$ and its essence is as follows.
First evaluate the speed of changing $Q(x(t))$ along trajectories of (\ref{sg1}), i.e. evaluate
$\dot Q=w(x,\theta,t)=\nabla Q(x)^{\trn}F\big(x, \theta,t \big)$. Next evaluate the gradient of this speed with respect to $\theta$: $\nabla_{\theta} w(x,\theta,t)$. Finally use the following algorithm
\begin{align}
\dot \theta =-\Gamma \nabla_{\theta} w(x,\theta,t),
\label{sg3}
\end{align}
where $\Gamma>0$ is a positive definite gain matrix.

Under some conditions, see \cite{F80,Fradkov90,FMN99,AF21} the algorithm (\ref{sg3}) ensures stabilization of the closed loop system in the following sense: the trajectories of system (\ref{sg1}), (\ref{sg3}) are bounded  and  goal (\ref{sg2}) is achieved.

\section{Robustification of the Speed Gradient algorithms with negative feedback}
In a variety of practical cases a realistic model of the controlled system dynamics should take into account the disturbances. In many cases the disturbances may be assumed bounded and act additively, i.e the plant model
instead of (\ref{sg1}) looks as follows
 \begin{align}
\dot x(t) = F\big(x, \theta,t \big) + f(t),
\label{sg4}
\end{align}
where $f(t)$ represent disturbances and $||f(t)||\le\Delta_f, \Delta_f>0.$
Influence of possibly nonvanishing disturbance in (\ref{sg4}) implies that it is not possible any more
to reach the goal (\ref{sg2}) and the goal should be relaxed. Since disturbance is bounded, the goal
may be reformulated as ultimate boundedness:
\begin{align}
\overline{\lim}_{t\to\infty}Q(x(t))\ge \Delta
\label{sg5}
\end{align}
for some $\Delta>0.$
The problem is to find an algorithm of the form
\begin{equation}\label{2.3}
	\dot{\theta}=\Theta(x,\theta,t)
\end{equation}
ensuring the goal (\ref{sg5}) for all $x(0),\theta(0)$ from some set of initial conditions $\Omega$.

Consider the  speed gradient adaptation algorithm with additive feedback in the form \cite{Fradkov90}:
\begin{equation}\label{2.43}
	\dot{\theta}=-\Gamma \left[\kappa \nabla_\theta w\left(x,\theta,t\right) + \zeta\left(\theta\right) \right],
\end{equation}
where $\zeta(\cdot): \mathbb{R}^m \rightarrow \mathbb{R}^m$-- is the feedback function, and the internal gain $\kappa>0$ is introduced for convenience of subsequent formulations.

In what follows we will also consider an algorithm of a more general form (combined),
\begin{equation}\label{2.44}
	\frac{d}{dt}\left[\theta+\psi\left(x,\theta,t\right)\right]=- \Gamma\left[k \nabla_\theta w\left(x,\theta,t\right) + \zeta\left(\theta+\psi\left(x,\theta,t\right)\right)\right]
\end{equation}
where $\psi(\cdot)$ satisfies pseudogradient condition
\begin{equation}\label{2.19}
	\psi(x,\theta,t)^{\trn}\nabla_{\theta}w(x,\theta,t)\ge 0.
\end{equation}
 The negativity of the parametric feedback $\zeta(\theta)$ means that it prevents excessive growth of the magnitude  $\|\theta(t)\|^2$ or (which is asymptotically the same -- of the magnitude $\|\theta(t)-\theta_*\|^2$, where $\theta_*$  is ``the ideal control''). Therefore the quantity $\zeta(\theta)^T (\theta-\theta_*)$ should be negative for sufficiently large  $\|\theta(t)-\theta_*\|$. Let's impose this requirement in a more convenient form.

Condition NF. There exist numbers $\rho>0,\rho'>0$ such that
$\zeta(\theta)^T(\theta-\theta_*)\geq \rho \|\theta-\theta_*\|^2 - \rho'. $

Condition NF is satisfied, e.g. for the function $\zeta(\theta)=\alpha\theta$, for which $\rho=\alpha/2$, $\rho'=\alpha\|\theta_*\|^2/2$ . One can also choose $\zeta(\theta)$  as follows $\zeta(\theta)=\alpha(\theta-\overline{\theta})$, where $\overline{\theta}\in\mathbb{R}^m$ is some prior estimate of the vector $\theta_*$. In this case $\rho'=\alpha\|\overline{\theta}-\theta_*\|^2/2$.

Before proceeding to formulate the robustness properties of algorithms (\ref{2.43}),(\ref{2.44}) let's discuss the requirements for setting the control objective under disturbance conditions. Firstly, when specifying the control objective (\ref{sg5}), the accuracy level $\Delta>0$ cannot be arbitrary but must depend on the level of disturbance $\Delta_\varphi>0$. Secondly, if it is known that some accuracy $\Delta_*>0 $ is achievable for system (\ref{sg1}) for some $\theta=\theta_*$, then for the system (\ref{sg4}), (\ref{2.3}) it is desirable to fix the goal with the same $(\Delta=\Delta_*)$ or with the arbitrarily close $(\Delta>\Delta_*)$ accuracy level. Thirdly it is desirable that the value  $\Delta_*$, be taken as the smallest possible value, e.g. $\Delta_*=J_*$; where $J_*=\inf\limits_{\theta}\sup\limits_{\varphi}\overline{\lim\limits_{t\rightarrow\infty}} Q(x(t),t)$;
this would give the system optimality or suboptimality properties. However, calculating $J_*$ in nonlinear problems is usually difficult. If one manages to calculate the value of
 $J_*$, it turns out to depend on parameters of the controlled system, which are unknown in adaptive control tasks due to the nature of the problem. The proposed way out of this situation is based on changing the formulation of the task: instead of approximating the exact optimum, one should consider approaching its estimate.

In accordance with above and following \cite{Fradkov90}, introduce key definitions.

{\bf Definition 1}. Consider the following class of the control problems $\Xi=\left\{\xi\right\}=\left\{F_\xi(\cdot),\varphi_\xi(\cdot),Q_\xi(\cdot),A_\xi\right\}$ and let for all $\xi \in \Xi$ there exist an algorithm $\Theta_\xi^* \in A_\xi$ such that
\begin{equation}\label{2.45}
J(\xi,\Theta_\xi^*) \leq \Delta_*,
\end{equation}
where $J(\xi,\Theta)=\overline{\lim\limits_{t\rightarrow\infty}} Q(x(t),t),x(t)$  is the solution of (\ref{sg4}), (\ref{2.3}). Then the value $\Delta_*$ is called the  estimate of the optimum for the functional $J(\xi,\Theta)$ in class $\Xi$. The control algorithm
  $\Theta_\xi^*$ , providing  (\ref{2.45}) for any $\xi$  
is called estimated  optimal in class $\Xi$. $\blacksquare$

{\bf Definition 2}. Let the class $\Xi$ of control problems,  the estimate $\Delta_*$ of the optimum of the functional $J(\xi,U)$  and the number $\xi>0$ be given. The control algorithm (\ref{2.3}), providing for any $\xi \in \Xi$ inequality $J(\xi,\Theta_\xi)\leq \Delta_* + \epsilon$ is called estimatedly $\epsilon$-optimal in class $\Xi$. $\blacksquare$

Note that in Definition 2 it is not required that the algorithm $\Theta_\epsilon$ would belong
to the set $A_\xi$, which, at first glance, makes finding an estimated  $\epsilon$-optimal algorithm trivial. In reality, additional requirements are imposed on $\Theta_\epsilon$.  For instance, in adaptive control problems, it is required that the sought algorithm does not depend on $\xi$,
i.e., it remains the same for all  $\xi \in \Xi$.

Additional constraints may be placed on measurable quantities, output controls, etc.

In order to formulate the main results  introduce a number of assumptions.

Assumption A1 (regularity). Functions $F(x,\theta,t),\nabla_xQ(x),\nabla_{\theta}w(x,\theta,t)$ are   bounded locally in $(x,\theta)$ uniformly in $t\ge 0$.

Assumption A2 (growth). Function $Q(x)$ is continuous, and grows unboundedly at infinity: $Q(x)\to+\infty$ as  $x\to\infty$.

Assumption A3 (convexity). Function $w(x,\theta,t)$ is convex in $\theta$, i.e. for
all $x,\theta,\theta',t, t\ge 0$ the inequality holds:
$$
w(x,\theta',t)-w(x,\theta,t)\ge(\theta'-\theta)^{\trn}nabla_{\theta}w(x,\theta,t).
$$

Assumption A4 (attainability). There exists vector $\theta_*$ and number $\alpha_0>0$ such that
$w(x,\theta_*,t)\le-\alpha_0Q(x)$.

Next consider the class  $\Xi$ of control problems defined by the following requirements:

- functions $F(\cdot),Q(\cdot)$ satisfy assumptions A1-A4 and in addition, the following quantitative assumptions hold.

Assumption A5. Function $Q(x)$ satisfies inequality $\|\nabla_xQ(x)\|\leq \alpha_1 Q(x)^\sigma$ for some $\alpha_1>0, 0\leq\sigma<1$;

Assumption A6. Function $f(t)$ is bounded: $\|f(t)\|\leq\Delta_f$;

Assumption A7. Set $A_\xi$ for all $\xi\in \Xi$ consists of algorithms producing constant controls: $\theta(t)\equiv \theta \in \mathbf{R}^m$.

In order to formulate the statement about properties of the robustified algorithm (\ref{2.43}) let us
formulate an auxiliary proposition.

{\bf Proposition 1}. Any number  $\Delta_*$ satisfying inequality
\begin{equation}\label{2.46}
	\left(\frac{\Delta_f\alpha_1}{\alpha_0}\right)^{1/(1-\sigma)} \leq \Delta_*
\end{equation}
is an estimate of the optimum of the functional
 $J(\xi, \Theta)=\overline{\lim\limits_{t\rightarrow\infty}} Q(x(t),t)$
 over the class $\Xi$.

To prove the proposition,  compute the rate of change of $Q(x(t), t)$ along system (\ref{sg4}), (\ref{2.3}) for $\Theta:\theta(t)\equiv\theta_*$, where $\theta_*\in \mathbf{R}^m$ is the ``ideal'' control from Assumption A4.

Using the easily verifiable inequality
$-\alpha_0 Q + \Delta_\varphi \alpha_1 Q^\sigma \leq -\alpha_0(1-\sigma)[Q-(\Delta_\varphi\alpha_1/\alpha_0)^{1/(1-\sigma)}]$,
yields
$\dot{Q_t} \leq \omega(x,\theta_*,t)+(\nabla_xQ)^Tf(t) \leq -\alpha_0Q + \Delta_f\alpha_1Q^\sigma \leq -\alpha_0(1-\sigma)(Q-\Delta_*)$
and therefore  $J(\xi,\Theta)\leq\Delta_*$.

Note that the inequality (\ref{2.46}) is exact: e.g. for
$n=m=1, F(x,\theta)=ax+b\theta$,
$Q(x)=|x|^\sigma , |f(t)| \leq \Delta_f$
inequality $J(\xi, \Theta)\leq (\Delta_f\alpha_1/\alpha_0)^{1/(1-\sigma)}$ becomes an equality
for $\theta(t)\equiv\Delta_f$

Accordingly, denote by  $\Xi_1(\Delta_*)$ the class of control problems
defined above, provided that (\ref{2.46}) holds for  given  $\Delta_*$.

The main property of SG-algorithm with negative feedback is that it ensures estimated $\epsilon$-optimality of system (\ref{sg4}), (\ref{2.44}) in the introduced class  $\Xi_1(\Delta_*)$ for any $\epsilon>0$.  The precise statement is formulated following \cite{Fradkov90} for the more general algorithm (\ref{2.44}). Then the statement for the algorithm (\ref{2.43}) will follow as a special case.

{\bf Theorem 1}.  For any problem from the class $\Xi_1(\Delta_*)$ using algorithm (\ref{2.44}) that satisfies conditions (\ref{2.19}) and NF, with $2\rho\lambda \geq \alpha_0(1-\sigma), \lambda=\lambda_{min}(\Gamma)$, the system (\ref{sg4}), (\ref{2.44}), is ultimately bounded. For $k>k_0$, where $k_0=2\rho'[\epsilon\alpha_0(1-\sigma)]^{-1}$, algorithm (\ref{2.44}) is {\it estimated $\epsilon$-optimal} in class $\Xi_1(\Delta_*)$.

{\it Corollary}. For any problem in $\Xi_1(\Delta_*)$ if the
conditions (\ref{2.19}), OS1, $2\rho\lambda \geq \alpha_0(1-\sigma)$ are fulfilled for the solutions of  system (\ref{sg4}), (\ref{2.44}) then the bound
$\overline{\lim\limits_{t\rightarrow\infty}} Q(x(t),t) \leq \Delta_* + 2 \rho' [k\alpha_0(1-\sigma)]^{-1}$ holds.

The drawback of the described robustification method is significant distortion of control performance for small $\nabla_\theta\omega(x,\theta,t)  ~ (\text{i.e. for} ~  \nabla_\theta\omega(\cdot)\equiv0, \zeta(\theta)=\theta- \overline{\theta}$ with algorithm (\ref{2.43}) $\lim\limits_{t\rightarrow\infty}\theta(t)=\overline{\theta}$ regardless of the history of process $x(t))$.

To overcome this shortcoming the feedback term $\zeta(\theta)$ is chosen to be zero on a bounded set containing the ``ideal'' control $\theta_*$.

Changing $\zeta(\theta)$ in a bounded set obviously does not affect the validity of condition NF, although the value of $\rho'$ may change. For example, one may choose
\begin{equation}
	\zeta(\theta) =
	\begin{cases}
		\alpha(\theta - \overline{\theta}), & \text{for } \|\theta - \overline{\theta}\| \geq d, \\
		0, & \text{for } \|\theta - \overline{\theta}\| < d.
	\end{cases}
\end{equation}

With this choice, the inequality in condition NF will hold for

$\rho=\alpha/2,~ \rho'=\alpha(~\|\overline{\theta}-\theta_*\|+d)^2/2.$

To improve the efficiency of the feedback law (2.47), it is advisable to increase the penalty for exiting the ball
 $\left\{\theta:\|\theta-\overline{\theta}\| \leq d\right\}$, i.e., to increase the value of  $\alpha$. Practical implementation considerations also lead us to introduce limits on  $\|\zeta(\theta\|)$ (saturation). By letting $\alpha\rightarrow\infty$, we arrive at a sign (relay) feedback 

\begin{equation}
		\zeta(\theta)=
		\begin{cases}
			 A~ sign(\theta-\overline{\theta}) & \text{for} \|\theta-\overline{\theta}\| \geq d, \\
		     0  & \text{for}  \|\theta-\overline{\theta}<d.
		     \end{cases}
\end{equation}

{\it Proof of Theorem 1}. Consider the function
$V(x,\theta,t)=Q(x,t)+0.5 k^{-1}\|\theta-\theta_*+\gamma\psi(x,\theta,t)\|^2\Gamma-1$
and evaluate its derivative with respect to the system (\ref{sg4}), (\ref{2.44}) :

$\dot{V}=w(x,\theta,t)+(\nabla_xQ)^Tf(t)+[\theta-\theta_*
+\gamma\psi(x,\theta,t)]^T \times \nabla_\theta w(x,,\theta,t) - k^{-1} [\theta-\theta_*+\gamma\psi(x,\theta,t)]^T \zeta(\theta+\gamma\psi(x,\theta,t)) \leq -\alpha_0Q +\alpha_1Q^\sigma \Delta_f- \frac{\rho\|\theta-\theta_*+\gamma\psi(x,\theta,t)\|^2}{k} + \frac{\rho'}{k}$.

Employing easily verifiable inequalities

$-\alpha_0Q +\Delta_\varphi \alpha_1 Q^\sigma \leq -\alpha_0(1-\sigma)(Q-\Delta_*)$,

$- \|\theta-\theta_*+\gamma\psi\|^2 \leq -\lambda\| \theta-\theta_*+\gamma\psi\|^2_{\Gamma^{-1}}$,

where $\lambda=\lambda_{min}(\Gamma)$, and taking into account
relation $2\rho\lambda \leq \alpha_0(1-\sigma)$, we arrive at linear differential inequality

$\dot{V}= \leq -\alpha V + \beta$, where
$\alpha=\alpha_0(1-\sigma,),   \beta=\beta(k)=\frac{\rho'}{k} + \alpha_0(1-\sigma)\Delta_*$.

Integration of the obtained differential inequality yields inequality

$Q(x(t),t) \leq V(x(t),\theta(t),t) \leq V(x(0),\theta(0),0) e^{-\alpha t} +\frac{\beta}{\alpha},$

implying $L$-dissipativity of the system. To prove estimated $\epsilon$-optimality
of the system let us fix $\epsilon>0$ and choose $k>k_0$. Then

$\beta(k) < \beta(k_0) = \alpha \left(\Delta_*+\frac{\epsilon}{2}\right)$
and therefore, the inequalities
$\overline{\lim\limits_{t\rightarrow\infty}} Q(x(t),t) \leq \Delta_* + \frac{\epsilon}{2} < \Delta_* + \epsilon$,
hold that prove the theorem.$\blacksquare$

In the end of this section it is worth to mention a different way of robustification based on introducing the deadzone into the adaptive control algorithm. The adaptation algorithm with the deadzone looks as follows:
\begin{equation}\label{deadzone}
		\dot\theta=
		\begin{cases}
			 -\Gamma\nabla_{\theta}w(x,\theta,t) & \text{for}  Q(x,t)\ge \Delta,\\
		     0  & \text{for}  Q(x,t)< \Delta.
		     \end{cases}
\end{equation}
The properties of the control algorithm with the deadzone (\ref{deadzone}) are formulated in the following theorem \cite{FFY}.

{\bf Theorem 2}. For any control problem from the class $\Xi{\Delta_*}$ all trajectories of the system
(\ref{sg4}), (\ref{deadzone}) are bounded and the goal (\ref{sg5}) is achieved for $\Delta_*<\Delta$.

\section{Example. Suboptimal robustness of adaptive passification based output feedback}

Consider a linear system affected by bounded disturbances

\begin{equation}\label{plant-dist}
\dot x=Ax+Bu+f(x,u,t), \quad y=L^Tx
\end{equation}
where $||f(x,u,t)||\le \Delta_f$. To solve the adaptive control problem
apply the passification SG approach \cite{F74,FFY,Fradkov90}. The designed
adaptive controller is as follows:
\begin{equation}\label{adap-out-con}
u=\theta^Ty, \dot \theta=-\Gamma[\kappa(g^Ty)y +\alpha(\theta-\bar\theta)],
\end{equation}
where $\kappa>0, \alpha>0$ and $\bar\theta\in R^l$ - apriori estimate for the vector of ``ideal'' parameters $\theta^*$. It follows from  Theorem 1 that
the system (\ref{plant-dist}),(\ref{adap-out-con}) is ultimately bounded ($L$-dissipative). If $\kappa>0$ is sufficiently large then the limit stabilization error is close to its minimum possible value (minimum is taken over all bounded disturbances such that $||f(x,u,t)||\le \Delta_f$).

Let $\Delta_x=\overline{\lim\limits_{t\rightarrow\infty}} ||x(t)||$, where $x(t)$ is a solution of (\ref{plant-dist}),(\ref{adap-out-con}) and let $Q(x)=0.5 x^THx$.
It is seen that function $Q(x)$ satisfies conditions of Theorem 1.
In this case the class $\Xi_1(\Delta_*)$ has the following  parameter values:
$\sigma=0.5, \alpha_1^2=\lambda^+, \alpha_0=2\sigma_*$, where $\lambda^+=\lambda_{max}(H), \sigma_*$ is the stability degree of matrix $A_*=A+B\theta_*^TL^T$. Evaluation of the accuracy bounds from Theorem 1 yields
$$
\Delta_x^2\le\overline{\lim}_{t\to\infty}\frac{Q(x(t))}{\lambda^{-}}\le\\
\frac{\Delta_f^2\alpha_1^2}{\alpha_0^2\lambda^{-}}+\frac{\epsilon}{\lambda^{-}},
$$
where $\lambda^{-}=\lambda_{min}(H), \epsilon>0$. Taking into account the simple inequality $\sqrt{p+q}\le\sqrt{p}+q/(2\sqrt{p})$, we obtain the final bound
\begin{equation}\label{error-bound}
\Delta_x\le \sqrt{\frac{\lambda^+}{\lambda^{-}}}\frac{\Delta_f}{2\sigma_*}+
\frac{\epsilon\sigma_*}{\Delta_f\sqrt{\lambda^+\lambda^-}}.
\end{equation}
Note also that the bound (\ref{error-bound}) holds if the following relations between the parameters of the adaptive controller are valid:
\begin{equation}\label{param}
\kappa>\frac{\alpha||\bar\theta-\theta^*||^2}{\epsilon\theta^*}, ~\alpha\lambda>0,
\end{equation}
where $\lambda=\lambda_{min}(\Gamma)$. Note that for the approximate
calculations the value of $\sigma^*$ can be taken as the stability degree
(minimal distance from the polynomial zeros to the imaginary axis) of the
polynomial $G(s)=g^Ta(s)$ which is the numerator of the transfer function $g^TW(p)$. It follows from the observation that $n-1$ eigenvalues of the matrix
$A^*=A+B\theta^{*T}L^{T}$ tend to zeros of $G(s)$, if $\theta^*=-\mu g,~\mu\to+\infty$.


\section{Conclusions}
This paper is dedicated to the problem of ensuring a suboptimal performance of adaptive systems under disturbances, i.e providing
performance arbitrarily close to that of an optimal system. The exposition is an updated version of the results in
 the book \cite{Fradkov90} where the concept of estimated optimality was introduced. It may be interesting for the readers since the book \cite{Fradkov90} is written in Russian
and a number of its results were never published in English. The goal of this note is to introduce the
updated results on estimated suboptimal adaptive control to the worldwide control community.

The first main result of the paper is evaluation the tight  accuracy bounds of speed-gradient control for a broad class of nonlinear nonlinearly parameterized  adaptive systems controlled by the speed-gradient algorithms
(Theorem 1 in Section 3). Then the case of adaptive control of a linear passifiable plant
by adaptive output feedback controller is studied as example.

Further research may be devoted to  numerical study of robustified algorithms for robotic systems
 and extending them to the nonsmooth case \cite{DolgFr21}.


\end{document}